\documentclass[12pt,reqno]{amsart}

\usepackage[final]{epsfig}
\usepackage{amssymb,amsmath,amstext,amsthm,amsfonts}
\usepackage{euscript,a4wide,pslatex}
\usepackage[ansinew]{inputenc}

\newtheorem{theorem}{Theorem}[section]
\newtheorem{proposition}[theorem]{Proposition}
\newtheorem{corollary}[theorem]{Corollary}
\newtheorem{lemma}[theorem]{Lemma}
\newtheorem{remark}[theorem]{Remark}
\newtheorem{maintheorem}{Theorem}

\numberwithin{equation}{section}

\newcommand{\bR}{{\bf R}}
\newcommand{\bS}{{\bf S}}

\newcommand{\bN}{{\bf N}}

\newcommand{\cA}{{\mathcal A}}
\newcommand{\cB}{{\mathcal B}}
\newcommand{\cC}{{\mathcal C}}

\newcommand{\cO}{{\mathcal O}}

\newcommand{\tofraco}{\rightharpoonup}

\newcommand{\de}{{\delta}}
\newcommand{\vfi}{{\varphi}}
\renewcommand{\epsilon}{\varepsilon}
\newcommand{\ep}{\varepsilon}

\newcommand{\te}{\theta}

\DeclareMathOperator{\inte}{int}
\DeclareMathOperator{\diam}{diam}

\DeclareMathOperator{\dist}{dist}
\DeclareMathOperator{\supp}{supp}

\theoremstyle{remark}

\title[Statistical stability of saddle-node arcs] {
Statistical stability of saddle-node arcs}

\begin{thanks} {V.A. was partially
    supported by CMUP-FCT (Portugal), post-doc grant
    SFRH/BPD/16082/2004/UX2U from FCT (Portugal) and CNPq
    (Brazil). Part of this work was done while enjoying a
    leave from CMUP at PUC-Rio/IMPA. M.J.P. was partially
    supported by CNPq-Brazil/Faperj-Brazil/Pronex Dyn.
    Systems.}
\end{thanks}

\author{V\'\i tor Ara\'ujo}
\address{Centro de Matem\'atica da
  Universidade do Porto, Rua do Campo Alegre 687, 4169-007
  Porto, Portugal}
\email{vdaraujo@fc.up.pt}
\urladdr{http://www.fc.up.pt/cmup/home/vdaraujo}

\author{Maria Jos\'e Pacifico}

\address{Instituto de Matematica, 
Universidade Federal do Rio de Janeiro, 
C. P. 68.530, CEP 21.945-970, 
Rio de Janeiro, R. J. , Brazil} 
\email{pacifico@impa.br and pacifico@im.ufrj.br}
\urladdr{http://www.dmm.im.ufrj.br}

\keywords{statistical and stochastic stability,
  equilibrium states, entropy formula, random perturbations,
  generic saddle-node unfolding}

\subjclass{Primary: 37A10, 37G25. Secondary: 37A35, 37H99.}

\begin{document}

\begin{abstract}
  We study the dynamics of generic unfoldings of saddle-node
  circle local diffeomorphisms from the measure theoretical
  point of view, obtaining statistical and stochastic
  stability results for deterministic and random
  perturbations in this kind of one-parameter families. 
\end{abstract}

\maketitle

\tableofcontents

\section{Introduction}
\label{sec:introduction}

The study of the modifications of the long term behavior
of a dynamical system undergoing perturbations of the parameters 
has been one of the main themes of Bifurcation Theory. In
the last decades the measure theoretical point of view has
been intensively developed emphasizing the understanding of
the asymptotic behavior of almost all orbits. The main notions
associated to this point of view are those of \emph{physical
  measure} and of \emph{stochastic or statistical stability}.

Let $M$ be a circle and $f_0:M\to M$ be a $C^2$ local
diffeomorphism.
An $f_0$-invariant probability measure $\mu$ is
 \emph{physical} if the \emph{ergodic basin}
\[
B(\mu)=\left\{x\in M:
\frac1n\sum_{j=0}^{n-1}\varphi(f_0^j(x))\to\int\varphi\,
d\mu\mbox{  for all continuous  } \varphi: M\to\bR\right\}
\]
has positive Lebesgue (length) measure in $M$. This means
that the asymptotic behavior of ``most points'' is
observable in a ``physical sense'' and determined by the
measure $\mu$.

Given a smooth family $(f_t)_{t\in[0,1]}$ of local
diffeomorphisms of $M$ admitting physical measures $\mu_t$
for every $t$, we say that $f_0$ is \emph{statistically
  stable} if $\mu_t$ tends to $\mu_0$ when $t\to0$ in a
suitable topology. This corresponds to stability of the long
term dynamics of most orbits under deterministic
perturbations of $f_0$. 

In this setting a straightforward consequence of the Ergodic
Theorem is that \emph{every ergodic $f_0$-invariant probability
measure $\mu_0$ absolutely continuous with respect to
Lebesgue measure $m$ is a physical measure}.

A \emph{random perturbation} of $f_0$ is defined by a family
of probability measures $(\te_\ep)_{\ep>0}$ on $[0,1]$ and the random
 sequence of maps
\[
f_{\omega}^n=f_{t_n}\circ\dots\circ f_{t_1},\quad
n\ge1\quad\mbox{and}\quad f^0_{\omega}=Id,
\]
where $Id:M\to M$ is the identity transformation, for
a sequence $\omega=(t_1,t_2,\dots)\in\supp(\te_\ep)^{\bN}$
and a given fixed $\ep>0$.
An invariant measure in this setting is called a
\emph{$\ep$-stationary measure}, which is a probability measure
$\mu$ such that for each continuous function $\vfi:M\to\bR$
\[
\int \vfi\, d\mu = \int \int \vfi(f_t(x))\, d\mu(x)\, d\te_\ep(t).
\]
Ergodicity in this setting needs an
extension of the notion of invariant set. We say that a
subset $E$ is \emph{$\ep$-invariant} when it satisfies
 \begin{itemize}
 \item[ ]if $x\in E$ then  $f_t(x)\in E$ for $\te_\ep$-almost every $t$,
 and
\item[ ] if $x\in M\setminus E$ then $f_t(x)\in M\setminus E$ for
  $\te_\ep$-almost every $t$.
 \end{itemize}
We say that a $\ep$-stationary measure $\mu$ is \emph{ergodic} if
$\mu(E)$ equals $0$ or $1$ for every $\ep$-invariant set $E$.
In this setting a point $x$ belongs to the \emph{ergodic
  basin} $B(\mu)$ if for all continuous $\varphi: M\to\bR$
and $\te_\ep^\bN$-almost every $\omega$ we have
\[
\frac1n\sum_{j=0}^{n-1}\varphi(f_\omega^j(x))\to\int\varphi\,
d\mu \quad \mbox{when} \quad n\to\infty.
\]
A stationary measure is \emph{physical} if the Lebesgue
measure of its ergodic basin is positive.
We again have that \emph{an absolutely continuous
ergodic stationary probability measure is physical}.

Assuming that the family $(\te_\ep)_{\ep>0}$ satisfies
$\supp(\te_\ep)\to\{0\}$ when $\ep\to0$ and there exist
physical stationary measures $\mu^\ep$ for every small
enough $\ep>0$, we say that $f_0$ is \emph{stochastically
  stable} if every limit point of $(\mu^\ep)_{\ep>0}$ when
$\ep\to0$ is a physical measure for $f_0$.  This corresponds
to stability of the asymptotic dynamics under random
perturbations of $f_0$.

In this paper we study the dynamics of generic unfoldings
$(f_t)_{t\in[0,1]}$ of a saddle-node circle local
diffeomorphism $f_0$ from the measure theoretical point of
view, obtaining statistical stability results for
deterministic and random perturbations in this kind of
one-parameter families. In particular we show that the map
is uniformly expanding for all parameters close enough to
the parameter of the saddle-node and has positive Lyapunov
exponent uniformly bounded away from zero.

This kind of results in the particular case of saddle-node
circle homeomorphisms might have applications to
the mathematical modeling of neuron
firing, see~\cite{pk02}.

Our results can be seen as an extension of the work
in~\cite{AT} where maps which are expanding everywhere
except at finitely many points were studied. Moreover these
results open the way into further study of the unfolding of
critical saddle-node circle maps considered in \cite{C03}.
In addition, piecewise smooth families unfolding a
saddle-node as in \cite{MPP} were used to build new kinds of
chaotic attractors for flows, and the statistical properties
of this kind of attractors can possibly be obtained through
suitable extensions of the techniques we present below.

\subsection{Statements of the results}
\label{sec:statements-results}

Let $f_0:M\to M$ be a $C^2$ local
diffeomorphism having a unique saddle-node fixed point that
we call $0$.

The fixed point $0$ is a \emph{saddle-node} if $f'(0)=1$ and
$f''(0)\neq0 (>0$ say).  A \textit{generic unfolding} of $0$
(or $f$) is a one-parameter family of maps $f_t:M\to M$ with
$t\in[0,t_0]$, so that $f_0=f$ and if $f(x,t)=f_t (x)$, then
$f(0,0)=0$, $\partial_x f(0,0)=1$, $\partial_x^2 f(0,0)>0$
and $\partial_t f(0,0)> 0$.  The family
$(f_t)_{t\in[0,t_0]}$ is called a \textit{saddle-node arc}
in \cite{DRV}.

Let $B(\{0\})$ be the basin of attraction of the saddle-node fixed
point $0$ for $f_0$, i.e.
\[
B(\{0\})=\{x\in M: f_0^k(x)\to 0 \quad\mbox{as}\quad k\to\infty \},
\]
and let the \emph{immediate basin} $W_0$ of $0$ be the
connected component of $B(\{0\})$ containing $0$.

\begin{figure}[htb] 
\begin{center}
\includegraphics[scale=0.6]{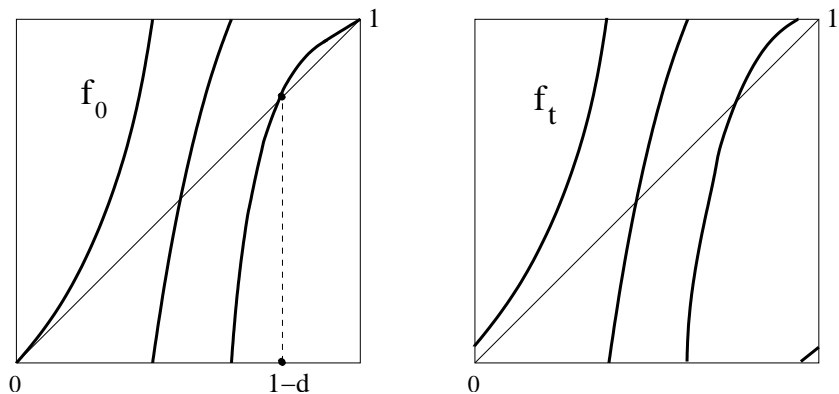}
\caption{\label{fig-vz0}
A saddle-node circle map.}
\end{center}
\end{figure}

We also assume the following global conditions on $f_0$,
\begin{description}
\item[H1] $|f'(x)|>1$ for all $x\in M\setminus W_0$,
\end{description}
see Figure~\ref{fig-vz0} for an example of such a map where
$W_0=[1-d,1]$. 

\begin{remark}
  \label{rmk:fonte}
We note that since $f_0$ is a local
diffeomorphism, there must be a fixed source $s$ ($s=1-d$ in
Figure~\ref{fig-vz0}) linked to the saddle-node, that is, a
connected component of $W^u(s)\setminus\{s\}$ is contained
in $W_0$.
\end{remark}

\begin{maintheorem}
  \label{thm:maindirac}
  Let $f_0$ be as above satisfying hypothesis (H1).  Then
  the Dirac mass $\de_0$ concentrated at $0$ is the unique
  physical measure of $f_0$.
\end{maintheorem}

The proof of this result in in
Section~\ref{sec:basin-attr-saddle}, where it is shown that
$B(\{0\})=M$ except for a zero Lebesgue measure subset of points.

\begin{maintheorem}
  \label{thm:main-FE}
  Let $f_0$ be as above satisfying hypothesis (H1).
  Then every $f_0$-invariant probability measure $\mu$
  satisfying the Entropy Formula
\begin{equation}
  \label{eq:E-F}
  h_{\mu}(f_0)=\int \log | f_0' | \, d\mu,
\end{equation}
must coincide with the Dirac mass $\de_0$ at the saddle-node
point $0$.
\end{maintheorem}

The proof of this theorem is in Section~\ref{sec:invar-meas-satisfy}.


\subsection{Statistical stability}
\label{sec:stat-stab-2}

The source linked to the saddle-node, see
Remark~\ref{rmk:fonte}, prevents the existence of either
sinks or nonhyperbolic period points in the unfolding of the
saddle-node. Using this we obtain the following statistical
stability result.

\begin{maintheorem}
  \label{thm:mainstatistical}
  Let $f_t:M\to M$ be a generic unfolding of $f_0$
  satisfying hypothesis (H1) above.    
  Then 
\begin{enumerate}
\item for every $t>0$ there exist $e_0=e_0(t)>0$ such that 
  \begin{enumerate}
  \item $f_t$ is uniformly expanding and there exists a
    unique absolutely continuous physical measure $\mu_t$
    whose basin equals $M$ except for a zero Lebesgue
    measure subset of points;
\item the Lyapunov exponent of Lebesgue almost every point
  is bigger  than $e_0(t)$.
  \end{enumerate}
\item $\mu_t\to\de_0$ when $t\to0$ in the weak$^*$ topology.
\end{enumerate}
\end{maintheorem}

We recall that item (2) means that \emph{$f_0$ is
  statistically stable} with respect to the unfolding given
  by $(f_t)_{t\ge0}$.

The proof of Theorem~\ref{thm:mainstatistical} is in
Section~\ref{sec:stat-stab}.

\subsection{Stability under random perturbations}
\label{sec:stab-under-rand}

Now we consider random perturbations of $f_0$ along the
family $f_t(x)=f_0(x)+t, \, x,t\in M$, which generically
unfolds the saddle-node at $0$, with a family
$(\te_\ep)_{\ep>0}$ of probability measures on $M$ such that
$\supp(\te_\ep)\to\{0\}$ when $\ep\to0$.

\begin{maintheorem}
  \label{thm:stochstabnosinks}
  Let $f_0$ satisfy hypothesis (H1) and let $f_t:M\to M$ be
  the family defined above. Then
  \begin{enumerate}
  \item for every family $(\te_\ep)_{\ep>0}$ as above
  satisfying additionally
  \begin{enumerate}
  \item $\te_\ep\ll m$;
  \item $\inte(\supp(\te_\ep))\neq\emptyset$;
  \item $\supp(\theta_\epsilon)\subset [0,t_0]$;
  \end{enumerate}
  for every $\ep>0$,  there exists a unique absolutely
  continuous stationary and ergodic probability $\mu^\ep$.
  \item $\mu^\ep\to\de_0$ when $\ep\to0$ in the weak$^*$ topology.
  \end{enumerate}
\end{maintheorem}

The above property (2) means that \emph{$f_0$ is
  stochastically stable} under absolutely continuous random
perturbations. 



\subsection{Statistical and stochastical stability for
  saddle-node circle homeomorphisms}
\label{sec.saddlehomeos}

Considering circle homeomorphisms with saddle-node points we
easily achieve the same results as we now explain.

\begin{figure}[htb] 
\begin{center}
\includegraphics[scale=0.6]{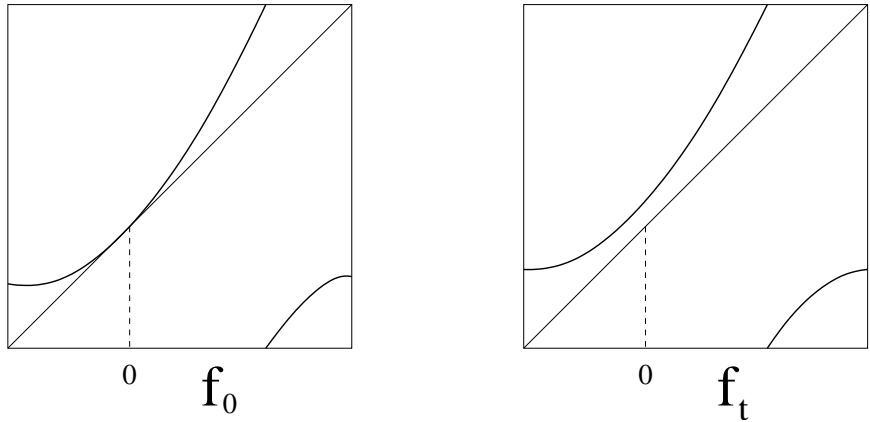}
\caption{\label{fig-vz1}
A saddle-node circle homeomorphism and a one-parameter family.}
\end{center}
\end{figure}

We say that a homeomorphism $f_0:M\to M$ is a 
\textit{saddle-node circle homeomorphism} if it satisfies
(see Figure \ref{fig-vz1}):
\begin{enumerate}
\item $f_0(0)=0$ and $f_0(x)\neq x$ for all $x\neq0$;
\item $f_0(x)>x$ for all $x\in V\setminus\{0\}$
for some open neighborhood $V$ of $0$.
\end{enumerate} 

Note that if $f_0$ were $C^2$ differentiable then
conditions  (1) and (2) above would imply that
$0$ was a usual $C^2$ saddle-node fixed point~\cite{NPT83}.

Since these kind of maps are \emph{uniquely ergodic} with measure
$\de_0$ (note that $f_0^n(x)\to 0$ when $n\to+\infty$ for
all $x\in M$) the following two stability results follow
from the fact that weak$^*$ accumulation measures of
stationary or invariant measures are invariant measures for
the limit map.

\begin{maintheorem}\label{circlet1}
  Let $f_0:M\to M$ be a saddle-node circle homeomorphism and
  let $(f_t)_{t\in[0,1]}$ be a continuous one-parameter
  family of circle homeomorphisms.  If we choose for every
  $t$ close to $0$ a $f_t$-invariant probability measure
  $\mu_t$, then every weak$^*$ accumulation point $\mu$ of
  the family $(\mu_t)_{t\in[0,1]}$, when $t\to 0$, is equal to
  the Dirac mass $\delta_0$ concentrated at the saddle-node.
\end{maintheorem}

This means that saddle-node circle homeomorphisms are
\textit{statistically stable}, i.e., the invariant
probability measure always vary continuously with the
unfolding parameter $t$ near the saddle-node parameter $0$.

Now let $(\te_\ep)_{\ep>0}$ be a family of probability
measures on $M$ for each $\ep>0$ such that
$\te_\ep\to\de_0$ when $\ep\to0$ in the weak$^*$ topology.

\begin{maintheorem}\label{circlet2}
  Let $f_0:M\to M$ be a saddle-node circle homeomorphism,
  $(f_t)_{t\in[0,1]}$ be a continuous one-parameter family
  of circle homeomorphisms and $(\te_\ep)_{\ep>0}$ be a
  family of probability measures on $M$ as above.

If we choose for every $\ep$ close to $0$ a stationary
probability measure $\mu^\ep$, then every weak$^*$
accumulation point $\mu$ of the family $(\mu^\ep)_{\ep>0}$,
when $\ep\to 0$, is equal to the Dirac mass $\delta_0$
concentrated at the saddle-node.
\end{maintheorem}

\bigskip

\noindent {\bf Acknowledgment.} We thank Y. Cao for pointing out to us
a mistake in the first version of this paper and for suggesting
a possible way to improve the statements of our results.


\section{Basins of attraction of sinks or  saddle-node points}
\label{sec:basin-attr-saddle}

Here we prove Theorem~\ref{thm:maindirac}. For this it is
enough to show that the basin $B(\{0\})$ of the saddle-node
$0$ has full Lebesgue measure in $M$.

\begin{theorem}
\label{th.baciatotal}
$m(M\setminus B(\{0\}))=0$.
\end{theorem}

In what follows  we set $g=f_0$  Clearly
to prove Theorem~\ref{th.baciatotal} it is sufficient to
obtain
\begin{proposition}
\label{pr.medida0}
\[
m\left(
I\cap
\bigcap_{n\ge0} g^{-n}(M\setminus W_0)
\right)
=0
\]
for every interval $I\subset M\setminus W_0$ whose length is small enough.
\end{proposition}
        
To prove this proposition we show that
for any given interval $I\subset M\setminus W_0$
there exists the first iterate $k$ such that 
$g^k(I)\not\subset M\setminus W_0$ and the relative
measure of the subinterval $G$ of points in $I$ that
fall into $W_0$ is a fixed proportion of the measure
of $I$.  For this we proceed as follows.

Let $I\subset M\setminus W_0$ be a given fixed interval and denote by
$\ell(I)$ its length. 
We observe that the boundary $\partial W_0$ of the immediate
basin consists of a source $s$.
This means that in a neighborhood outside $W_0$ we
always have some expansion.

For $\eta>0$ small we define the
following compact subset
\[
W(\eta)=\{ x\in W_0: d(x,M\setminus W_0)\ge\eta\}.
\]
We assume that $\ell(I)\le 1/4$ (recall that $\ell(M)=1$). Let us
choose $\eta_0>0$ small enough such that
\begin{equation}
\label{eq:1}
\int_J |g'| > 1
\end{equation}
for every interval $J\subset M\setminus W_0$ whose length
equals $\ell(I)$ and with one endpoint in 
$\partial W(\eta_0)$.
Then there exists $\sigma>1$ such that
\begin{equation}
\label{eq:2}
\int_J |g'| \ge\sigma
\end{equation}
for every interval $J\subset M\setminus W(\eta_0)$ such that
$\ell(J)\ge\ell(I)$.

\begin{remark}
\label{rmk.sigma-r}
The value of $\sigma$ depends on $\eta$ but if $0<\eta<\eta_0$ then
$\sigma(\eta_0)=\sigma(\eta)$.
\end{remark}

This uniform rate of expansion ensures the following.
\begin{lemma}
\label{le.1}
For any $0<\eta<\eta_0$ there exists $k_1$ such that
\[
g^k(I)\subset M\setminus W(\eta), \quad k=0,\dots,k_1-1
\quad\mbox{and}\quad 
g^{k_1}(I)\not\subset M\setminus W(\eta).
\]
\end{lemma}

\begin{proof}
We define
\[
L_0=\max\{ \ell(C): C \mbox{ is a connected component of }
M\setminus W(\eta)\}.
\]
If
$g^k(I)\subset M\setminus W(\eta), \quad k=0,\dots,k_0-1$ for some
$k_0>0$, we obtain
$\ell(g^k(I))\ge\sigma\ell(g^{k-1}(I))$ for all $1\le k\le k_0$.
Thus $\ell(g^{k_0}(I))\ge\sigma^{k_0}\ell(I)$.

If $k_0$ were arbitrarily large, then we would have
\begin{equation}
\label{eq.casoL}
\ell(g^{k_0}(I))\ge\sigma^{k_0}\ell(I)>L_0.
\end{equation} 
Thus by definition of $L_0$ we must have
 $g^{k_0}(I)\not\subset M\setminus W(\eta)$ as stated.
\end{proof}

Now it is easy to see that after a finite number of iterates
either $g^{k_1}(I)$ is completely inside the basin of the saddle-node,
or it contains a piece of the basin of uniform size $\eta$.

\begin{lemma}
\label{le.conexo}
If $k_1$ is given by Lemma~\ref{le.1} then
\begin{itemize}
\item either $g^{k_1}(I)\subset W_0$
\item or $g^{k_1}(I)\cap (M\setminus W_0)\neq\emptyset$ and
$g^{k_1}(I)\cap W(\eta)\neq\emptyset$.
\end{itemize}
Moreover in the last case we have $\ell(g^{k_1}(I)\cap
(W\setminus W(\eta)))\ge
\eta$.
\end{lemma}

\begin{proof}
The lemma follows from the fact that $g^{k_1}(I)$ is connected.
\end{proof}

Now we use the following bounded distortion result to estimate
the size of the piece of $I$ which is sent into $W_0\setminus W(\eta)$.

\begin{lemma}[Bounded distortion]
\label{le.boundist}
For $I\subset M\setminus W(\eta)$, $k_1$ given by
Lemma~\ref{le.1}
and $x,y\in I$ it holds
\begin{equation}
\label{eq.0}
\log\left| \frac{(g^{k_1})'(x)}{(g^{k_1})'(y)}\right|
\le C_0
\quad\mbox{where}\quad
C_0=\sup\left|\frac{g''}{g'}\right|\cdot
\frac1{1-\sigma^{-1}}.
\end{equation}
\end{lemma}

\begin{proof}
  Since $g$ is a local diffeomorphism , if $g^k(I)\subset
  M\setminus W(\eta), \quad k=0,\dots,k_1-1$ for some
  $k_1>0$ given by Lemma~\ref{le.1}, then by the definition
  of $\sigma$ we get $\ell(g^k(I))\ge\sigma\ell(g^{k-1}(I))$
  for all $1\le k\le k_1$.  We have
\begin{eqnarray*}
\log\left| \frac{(g^{k_1})'(x)}{(g^{k_1})'(y)}\right|
&=&
\sum_{j=0}^{k_1-1}
\left|
\log g'(g^j(x)) - \log g'(g^j(y))
\right|
=
\sum_{j=0}^{k_1-1} | (\log g')'(z_j) | \,
\ell([g^j(x),g^j(y)])
\\
&\le&
\sup\left|\frac{g''}{g'}\right|
\sum_{j=0}^{k_1-1} \ell(g^j(I))
\le
\sup\left|\frac{g''}{g'}\right|
\sum_{j=0}^{k_1-1} \sigma^{-(k_1-j)}\ell(g^{k_1}(I))
\\
&\le&
\sup\left|\frac{g''}{g'}\right|\cdot
\ell(M)\cdot
\frac1{1-\sigma^{-1}}
=\sup\left|\frac{g''}{g'}\right|\cdot
\frac1{1-\sigma^{-1}}.
\end{eqnarray*}
\end{proof}

\begin{corollary}
\label{co.1}
There exists a constant $C>0$ such that for every
interval $G\subset I$ and for $k_1$ given by Lemma~\ref{le.1}
we have
\[
\frac1{C}\cdot\frac{m(G)}{m(I)}
\le
\frac{m(g^{k_1}(G))}{m(g^{k_1}(I))}
\le C \cdot  \frac{m(G)}{m(I)}.
\]
\end{corollary}

\begin{proof}
It is straightforward to write for some $z\in I$
\begin{eqnarray*}
\frac{m(g^{k_1}(G))}{m(g^{k_1}(I))}
&=&
\frac{\int_G [(g^{k_1})'(x) / (g^{k_1})'(z)] \,dm(x) }
{ \int_I [(g^{k_1})'(x) / (g^{k_1})'(z)] \,dm(x) }
\le
C_0^2 \cdot \frac{m(G)}{m(I)}.
\end{eqnarray*}
Analogously we get 
\[
\frac{m(g^{k_1}(G))}{m(g^{k_1}(I))}\ge
\frac1{C_0^2} \cdot \frac{m(G)}{m(I)}
\]
showing that the corollary holds with $C=C_0^2$.
\end{proof}

Now we are ready to exclude from $I$ the points that
fall into the basin of the saddle-node in a controlled
way.

Let $G=(g^{k_1}\mid I)^{-1} (W_0)\subseteq I$. 
On the one hand, if $G\neq I$ then by
Lemma~\ref{le.conexo} and Corollary~\ref{co.1} we obtain
\begin{equation}
\label{eq.6}
\frac{m(G)}{m(I)}
\ge
\frac1{C}\frac{m(g^{k_1}(G))}{m(g^{k_1}(I))}
\ge
\frac{m(g^{k_1}(I)\cap(W_0\setminus W(\eta)))}
{C\cdot \sup|g'| \cdot m(M\setminus W(\eta))}
\ge
\frac{\eta}{C\cdot \sup|g'| \cdot m(M\setminus W(\eta))}
\end{equation}
where by definition of $k_1$ we have
\[
m(g^{k_1}(I))\le \sup|g'|\cdot m(g^{k_1-1}(I))
\le \sup|g'| \cdot m(M\setminus W(\eta)).
\]
Taking $\eta>0$ small enough (see also
Remark~\ref{rmk.sigma-r})~\eqref{eq.6} gives
\begin{eqnarray*}
m(I\setminus G) & = & m(I)-m(G) \le
m(I) \left(
1- \frac{\eta}{ C \cdot \sup|g'| \cdot m(M\setminus W(\eta))}
\right)
\\
&\le&
\gamma_0 \cdot m(I),
\end{eqnarray*}
where $\gamma_0\in (0,1)$ does not depend on $I$ nor on
$k_1$.

On the other hand, if $G=I$ then the last inequality is
trivially true and we are done.

\medskip

Otherwise, if a positive Lebesgue measure set remains in
$I\setminus G$, we proceed by induction to conclude the
proof of Proposition~\ref{pr.medida0}.

In what follows we set $I_0=I$ and 
$I_1=I\setminus G$.
Let us assume that we have already constructed a nested collection
of sets $I_0\supset I_1\supset \dots\supset I_n$ such that
\begin{enumerate}
\item for each
$i=1,\dots, n$, $I_i$ is a collection of intervals $J_{i,j}$
contained in $I_{i-1}$ and
\item to each $J_{i-1,j}$ there corresponds
an integer $k=k(i-1,j)\in\bN$ and a value $\eta=\eta(i-1,j)\in(0,\eta_0)$
satisfying (recall \eqref{eq:1})
\[
g^l(J_{i-1,j})\subset M\setminus W(\eta)
\mbox{  for all  }
l=0,\dots,k-1
\mbox{  and  }
g^k( J_{i-1,j}\setminus I_i)\subset W_0.
\]
\end{enumerate} 

The previous lemmas show that the following
result is true.

\begin{lemma}
\label{le.welldefined}
The sequence $I_n$ is well defined for all $n\ge1$
(it can be empty from some value of $n$ onward)
and
\[
m(I_{n+1})\le \gamma_0 \cdot m(I_n).
\]
\end{lemma}

We conclude that $m(\cap_{n\ge0} I_n)=0$. We now
show that this implies Proposition~\ref{pr.medida0}.

\medskip

Let us take $x\in\cap_{n\ge0} I_n$.
Then there exists a sequence $0=k_0<k_1<k_2<\dots$ of
integers and $\eta_1,\eta_2,\dots$ of reals in $(0,\eta_0)$ such that
\[
g^j(x)\in M\setminus W(\eta_i)
\mbox{  for  }
k_i\le j < k_{i+1}
\mbox{  and  }
i\ge0.
\]
Moreover $M\setminus W(\eta_i)\supset M\setminus W(\eta_0)$
for all $i\ge0$.  Hence $x\in g^{-j}(M\setminus W(\eta_i))
\supset g^{-j}(M\setminus W(\eta_0))$.

We deduce that if $g^j(y)\in M\setminus W(\eta_0)$ and $y\in I$, then
$y\in g^{-j}(M\setminus W(\eta_0))\subset g^{-j}(M\setminus W(\eta_i))$
and thus $y\in\cap_{n\ge0} I_n$. This means that
\[
I\cap
\bigcap_{n\ge0} g^{-n}(M\setminus W)
\subset
I\cap\left(
\bigcap_{n\ge0} g^{-n}\big( M\setminus W(\eta_0) \big)
\right)
\subset \bigcap_{n\ge0} I_n.
\]
Since we already know that $m(\cap_{n\ge0} I_n)=0$,
this ends the proof of Proposition~\ref{pr.medida0}.



\section{Invariant measures satisfying the Entropy Formula}
\label{sec:invar-meas-satisfy}

Here we prove Theorem~\ref{thm:main-FE}.  Let $\mu_0$ be a
$f_0$-invariant probability measure satisfying the Entropy
Formula~\eqref{eq:E-F}, i.e., $\mu_0$ is a equilibrium state
for the potential $-\log|f_0'|$. The following result shows
that we can assume without loss that $\mu_0$ is ergodic.

\begin{lemma}\label{lem.ergdecompequil}
 Almost every ergodic
component of an equilibrium state for $-\log|f_0'|$ is
itself an equilibrium state for this same function.
\end{lemma}

\begin{proof}
Let $\mu$ be an $f_0$-invariant measure satisfying
$h_\mu(f_0)=\int \log|f_0'| \, d\mu$. On the one hand,
the Ergodic Decomposition Theorem (see e.g
Ma\~n\'e~\cite{Man87}) ensures that
\begin{equation}
  \label{eq:9}
\int \log|f_0'| \, d\mu = 
\int\!\! \int \log|f_0'| \, d\mu_z \, d\mu(z) 
\quad\mbox{and}\quad h_\mu(f_0)=\int
h_{\mu_z}(f_0) \, d\mu(z).
\end{equation}
On the other hand, Ruelle's inequality guarantees for a
$\mu$-generic $z$ that 
\begin{equation}
  \label{eq:10}
h_{\mu_z}(f_0)\le\int \log|f_0'| \, d\mu_z.
\end{equation}
By~\eqref{eq:9} and~\eqref{eq:10}, and because $\mu$ is an
equilibrium state, we conclude that we have
equality in~\eqref{eq:10} for $\mu$-almost every $z$.
\end{proof}

Now we have two cases.
\begin{enumerate}
\item If $\mu_0(\{0\})>0$ then $\mu_0=\de_0$ because $\mu_0$ is
  ergodic and $0$ is fixed.
\item Else if $\mu_0(\{0\})=0$ then we let $x$ be a $\mu_0$-generic
  point, that is
\[
\frac1n\sum_{j=0}^{n-1} \de_{f_0^j(x)} \tofraco \mu_0
\quad\mbox{when}\quad n\to\infty,
\]
and we subdivide the argument in two more cases.
\begin{enumerate}
\item Either $x\in B(\{0\})$ or
\item $x\notin B(\{0\})$.
\end{enumerate}
\end{enumerate}
In case (a) since $x$ is a $\mu_0$-generic point we conclude that
$\mu_0=\de_0$ also.

In case (b) the positive orbit $\cO^+_{f_0}(x)$ of $x$ is contained in
the region of $M$ where $|f_0'|>1$, thus the integral in the Entropy
Formula is positive and so $h_{\mu_0}(f_0)>0$. 

It is known \cite{QZ01} that measures satisfying the Entropy Formula
with positive entropy for endomorphisms of one-dimensional manifolds
must be absolutely continuous (with respect to Lebesgue (length)
measure). 

Finally, since by Theorem~\ref{th.baciatotal} we have
$B(\{0\})=M, m\bmod0$, the absolute continuity of $\mu_0$
implies that there exits a $\mu_0$-generic point $x$ in
$B(\{0\})$, thus $\mu_0=\de_0$ as we wanted, proving
Theorem~\ref{thm:main-FE}.


\section{Statistical stability}
\label{sec:stat-stab}

Here we prove Theorem~\ref{thm:mainstatistical}.
First we recall some properties of the generic unfolding of
saddle-node arcs, which can be found in~\cite{DRV,NPT83}.

\subsection{Transition maps for saddle-node unfoldings}
\label{sec:trans-maps-saddle}

In what follows we let $f_0$ be a saddle-node local
diffeomorphism and perform a local analysis of the dynamics
near the saddle-node point $0$.  In this
setting the map $f_0$ is a $C^2$ diffeomorphism in a
neighborhood of $0$.

Given a saddle-node arc $(f_t)_t$ of one dimensional maps,
as defined in Section~\ref{sec:statements-results}, there is
what is called \emph{an adapted arc of saddle-node vector fields}
$(X(t,.))_t$, which has the form
\begin{equation}
\label{normalform}
X(t,x)=t+\alpha x^2+\beta x t+\gamma t^2+O(|t|^3+|x|^3),
\,\,\mbox{with $\alpha >0$},
\end{equation}
and describes the local dynamics of $(f_t)_t$:
the arc $(f_t)_t$ embeds as the time-one of $(X(t,.))_t$.
That is, if $X_s(t,.)$ denotes the time-$s$ map induced by 
$(X(t,.))_t$ then
$
f_t (x)=X_1(t , x)
$
for every $t$ and every $x$.
For $a<0<b$ fixed close enough to $0$, $k\in\bN$
and $t>0$ sufficiently small, if $\sigma_k(t)\in [0,1]$
is defined by the relation
$$
X_{k+\sigma_k(t)}(t,a)=b,
$$
then it is proved in \cite{DRV} that  for $k\ge 1$ large enough, 
there is a unique $t_k^*>0$ such that $\sigma_k(t_k^*)=0$, and 
$$
\sigma_k:[t_{k+1}^*,t_k^*]\rightarrow [0,1]
$$
is a $C^{\infty}$ decreasing diffeomorphism onto $[0,1]$. 
Set $t_k$ the inverse of $\sigma_k$.

For each $k\ge 1$ large enough, define 
$T_k:[0,1]\times [f_0^{-1}(a),f_0(a)]\to\bR$ by
$
T_k(\sigma,x)=f_{t_k(\sigma)}^k(x).
$
Note that $T_k$ depends on both $a$ and $b$.
For $ f_0^{-1}(a) < x < f_0(a)$ and $t$ small, define 
$t_a(t,x)$ by $X_{t_a(t,x)}(t,x)=a.$
The sequence $(T_k)_k$ converges in the $C^\infty$ 
topology to the {\em transition map} 
$$T_\infty : [0,1] \times [f_0^{-1}(a),f_0(a)]\mapsto\bR,
$$ 
defined by
$
T_\infty(\sigma,x)= X_{t_a(0,x)-\sigma}(0,b).
$
Note that $T_\infty$ depends also on both $a$ and $b$.

Observe also that $\partial_x T_\infty(\sigma,x)$ is bounded away 
from zero by a constant which does not depend on $(\sigma,x)$.
With $b$ fixed,  and taking  $a$ sufficiently close to $0$, 
we can assume that this constant is arbitrarily large, since
the number of iterates needed to take $a$ to $b$ increases
without limit if $t$ is small enough and $a$ close enough to $0$.


\subsection{Uniform expansion}
\label{sec:stat-stab-1}

Now we present the arguments proving statistical stability
of saddle-node arcs.

As in the previous subsection we fix $a<0<b$ with $a$ close
enough to $0$ in order to get $\partial_x T_a(\sigma,x)\ge
2 c_0>1$, for every $\sigma\in [0,1]$ and for all
$x\in[f_0^{-1}(a),f_0(a)]$. For small $t>0$ there exists
$k\ge1$ such that $t\in[t_{k-1}^*,t_k^*]$ and
\[
T_{k+\sigma_k(t)}: [f_0^{-1}(a),f_0(a)] \to
[f_0^{-2}(b),\infty), \quad x\mapsto f^k_{t_k(t)}(x)
\]
has derivative bigger than $c_0>1$, i.e.
\begin{equation}
  \label{eq:38}
  (T_{k+\sigma_k(t)})' \ge c_0 >1.
\end{equation}

\begin{remark}
  \label{rmk:antes4.1}
We have $f_t(x)=f(x)+t>f(x)$ for all $t>0$ and $x$ in a
small neighborhood of $0$. Since $f_0''(0)>0$ this ensures
that if $x$ is near $0$ and
$\omega_1,\dots,\omega_k\in[t_0^-,t_0]$ with $t_0^->0$, then
\[
\big( f_\omega^k\big)'(x)=\prod_{i=1}^k
f_0'\big(f_{\omega_i}\circ\dots\circ f_{\omega_1}(x)\big)
\ge \big( f_{t_0^-}^k\big)'(x).
\]
Hence the derivative of the transition maps
$T_{k+\sigma_k(t)}$ can be used as a lower bound for the
derivative along random orbits near $0$.
\end{remark}

\begin{theorem}
\label{thm.expansao}
There exist $t_0>0$ small enough such that for every
probability measure $\theta$ supported in $[0,t_0]$ with
$\theta(\{0\})<1$ and for every $x\in \bS^1$ 
\[
\limsup_{n\to\infty} \frac1n \log\Big| \big( f_\omega^n
\big)'(x) \Big| \ge0\quad \mbox{for}\quad
\theta^{\bN}-\mbox{a.e.  } \omega\in[0,t_0]^{\bN}.
\]
Moreover for every $0<t_0^-<t_0$ there exists
$e_0=e_0(t_0^-,t_0)>0$ such that for all $\omega
\in[t_0^-,t_0]^{\bN}$ and every $x\in\bS^1$
\begin{equation}
\label{eq.liap}
 \limsup_{n\to+\infty} \frac1n \log
 |(f_\omega^n)'(x)| \ge e_0.
\end{equation}
\end{theorem}

\begin{proof}
  To obtain such result we note that since we are assuming
  that $f_t$ is expanding outside the immediate basin of the
  saddle-node, it is enough to analyse the dynamical
  behavior near $0$. 

We fix $d_0\in W_0$ very close to the source $s$ connected
to the saddle-node, where $W_0$ is the immeadiate basin of
attraction of the saddle-node $0$ for $f_0$ (see
Section~\ref{sec:basin-attr-saddle} and
Remark~\ref{rmk:fonte}). We note that for all $t\in\bS^1$ we
have
\begin{equation}
  \label{eq:expand0}
  \Big| f_t'(x) \Big| \ge \sigma_0\quad\mbox{for all}\quad
x\in\bS^1\setminus[d_0,f_0^{-1}(a)],
\end{equation}
for some $\sigma_0>1$. Now we fix $x\in\bS^1$ and define for
any given $\omega\in[0,t_0]^\bN$
\[
R=R(\omega,x)=\{ i\ge0 : f_\omega^i(x)\in[d_0,f_0^{-1}(a)]\}.
\]
If $R=\emptyset$, then by~\eqref{eq:expand0} we have $\Big|
\big(f_{\omega_i}\big)'(f_\omega^i(x))\Big|\ge\sigma_0$ for
every $i\ge1$ and thus \eqref{eq.liap} holds with
$e_0=\log\sigma_0$. Otherwise $R\neq\emptyset$ and we set
$k=\min R$.

If $k>0$, then
$\Big|\big(f_{\omega_i}^k\big)'(x)\Big|\ge\sigma_0^k$ by
construction. Otherwise $k=0$ and so $x\in[d_0,f_0^{-1}(a)]$.
In this case we set
\[
\ell=\ell(\omega,x)=
\min\{i>0: f_\omega^i(x)\in[f_0^{-1}(b),f_0(b)]\}.
\]
We consider first the case $\omega\in[t_0^-,t_0]^\bN$ with
$0<t_0^-<t_0$. In this setting there is
$\ell_0=\ell_0(t_0^-)$ such that $\ell(\omega,x)\le\ell_0$
for all $x\in[d_0,f_0^{-1}(a)]$. We note that
$\ell_0(t_0^-)\to+\infty$ when $t_0^-\to0^+$.

The transition map and the geometry near the saddle-node
(see Remark~\ref{rmk:antes4.1}) ensure that there is
$c_\ell>1$ such that
\[
\Big|\big(f_{\omega}^\ell\big)'(x)\Big|\ge c_\ell>1\quad
\mbox{for}\quad
x\in[d_0,f_0^{-1}(a)], \, \omega\in[t_0^{-},t_0]^\bN, \,
\ell=\ell(\omega,x).
\]
If we define
\[
\beta=\beta(\ell_0)=
\min\Big\{\frac1\ell\log c_\ell : \ell=1,\dots,\ell_0 \Big\},
\]
then $\Big| \big(f_\omega^\ell\big)'(x)\Big| \ge
(e^\beta)^\ell$ and we set $\sigma_1=e^\beta$.

We have shown that for all
$(\omega,x)\in[t_0^-,t_0]^\bN\times\bS^1$ there is a
sequence $n_1<n_2<n_3<\dots$ satisfying
\begin{equation}
  \label{eq:mais1}
  \Big| \big(f_\omega^{s_k}\big)'(x)\Big| \ge \sigma^{s_k}
  \quad\mbox{with}\quad
s_k=n_1+\dots+n_k\quad\mbox{and}\quad
  \sigma=\min\{\sigma_0,\sigma_1\}. 
\end{equation}
Here $\sigma$ depends on $t$ through $\ell_0$ and
$\beta$. Hence~\eqref{eq.liap} holds with $e_0=\sigma$,
finishing the proof of the second part of the statement.

For the first part of the statement, we may assume that a
$\theta^\bN$-generic $\omega$ does not contain an infinite
sequence of coordinates arbitrarily near $0$, which is
enough to deduce that $\ell(\omega,x)$ is finite for every
$x\in\bS^1$ and $\theta^\bN$-a.e. $\omega$. Let us be more
precise.

The assumption $\theta(\{0\})<1$ ensures that for every
$\ep>0$ there exists $0<\delta_0<\epsilon$ such that
$\theta([0,\delta])<1$ for all $0\le\delta<\delta_0$. Hence
$X_0=[0,\delta]^\bN$ and $X_n=[0,t_0]^n\times X_0, n\ge1$ satisfy
$\theta^\bN(X_n)=0$ for all $n\ge0$ and thus
$Y=[0,t_0]^\bN\setminus \cup_{n\ge0} X_n$ has full
$\theta^\bN$-measure. In particular, letting $\delta=0$ we
get that a $\theta^\bN$-generic $\omega$ has no zeroes.

This shows that a $\theta^\bN$-generic sequence $\omega$
admits $\delta>0$ and a subsequence $n_1<n_2<n_3<\dots$ such
that
\[
\omega_{n_k}>\delta\quad\mbox{for all}\quad
k\ge1\quad\mbox{and}\quad \omega_n>0\quad\mbox{for all}\quad n\ge1.
\]
We conclude that $\ell(\omega,x)<+\infty$ (although it may
be arbitrarily big) for all $x\in\bS^1$ and
$\theta^\bN$-a.e. $\omega$.

Finally, going back to the initial argument, we assume that
$R\neq\emptyset$ and $k=0$, and set
\[
\alpha= \min\{ \omega_1,\dots,\omega_{\ell(\omega,x)}\} >0.
\]
Again by Remark~\ref{rmk:antes4.1} we see that there is
$c=c(\alpha)$ such that
\[
\Big| \big(f_\omega^{\ell(\omega,x)}\big)'(x)  \Big| \ge c
>1
\quad\mbox{for}\quad x\in[d_0,f_0^{-1}(a)],
\omega\in[\alpha,t_0]^\bN. 
\]
Since $\ell(\omega,x)$ can be arbitrarily big, the exponent
$\ell(\omega,x)^{-1}\cdot \log c$ can be arbitrarily close
to zero. Therefore the value of $\sigma$ in~\eqref{eq:mais1}
must be $1$, finishing the proof of the theorem.

\end{proof}

\begin{remark}
  \label{rmk:mais1}
The proof of Theorem~\ref{thm.expansao} shows, in
particular, (discarding the iterates near $0$) that for all
$x\in\bS^1$ and $\theta^\bN$-a.e. $\omega\in[0,t_0]^\bN$
there exists a sequence $n_1<n_2<n_3<\dots$ such that for
all $k\ge1$
\begin{description}
\item[a] $\Big|\big(f_\omega^{n_{2k-1}}\big)'(x)\Big|\ge1$;
  and
\item[b]
  $\Big|\big(f_{\sigma^{n_{2k-1}}\omega}^{n_{2k}-n_{2k-1}}\big)'\big(
  f_\omega^{n_{2k-1}}(x) \big)\Big|\ge\sigma_0^{n_{2k}-n_{2k-1}}$.
\end{description}
\end{remark}

Since Theorem~\ref{thm.expansao} ensures, in particular,
that for every $x\in M$ and $t\in(0,t_0)$ there exists
$n(t)\ge1$ such that $|(f_t^{n(t)})'(x)|>1$, then we
conclude that $f_t$ is uniformly expanding, i.e., there are
constants $C(t)>0$ and $\sigma(t)>1$ satisfying
$|(f_t^k)'(x)|\ge C(t)\cdot\sigma(t)^k$ for all $x\in M$, $k\ge1$
and every given $t\in(0,t_0)$, see e.g. \cite{AAS}.

\begin{theorem}
\label{thm.expansao1}
Let $t_0>0$ be given by Theorem~\ref{thm.expansao}.
Then for all $t\in(0,t_0)$ there exists a unique absolutely
continuous ergodic probability measure $\mu_t$ for $f_t$
such that
\begin{equation}
  \label{eq:entropyformula}
0< h_{\mu_t}(f_t)=\int \log | f_t' | \, d\mu_t.  
\end{equation}
\end{theorem}

\begin{proof}
The conclusion of Theorem~\ref{thm.expansao} is enough to
guarantee that $f_t$ is uniformly expanding, for each fixed
$t\in(0,t_0)$, by~\cite[Theorem A]{AAS}. It is well known
that uniformly expanding maps admit a unique absolutely
continuous ergodic invariant measure satisfying the Entropy
Formula~\eqref{eq:entropyformula}, see e.g.~\cite{Man87}.
\end{proof}

Another consequence of uniform expansion is the following.

\begin{theorem}
\label{thm.particao}
Let $\mu_0$ be a weak$^*$ accumulation point of $\mu_t$ when $t\to0$.
Then 
there exists a finite partition $\xi$ of $M$ which is
a $\mu_t\bmod0$
generating partition for $f_t$, for all $t\in(0,t_0)$, and
also that $\mu_0(\partial\xi)=0$, i.e., the $\mu_0$ measure
of the boundary of the atoms of $\xi$ is zero.
\end{theorem}

\begin{proof}
  Any finite partition of $M$ Lebesgue modulo zero is a
  $\mu_t\bmod0$ partition of $M$ (since $\mu_t\ll m$) and
  also a generating partition, by the uniform expansion of
  $f_t$ for $t\in(0,t_0)$, see e.g.~\cite{Man87}.
  
  A finite partition Lebesgue modulo zero whose boundary has
  also zero measure with respect to $\mu_0$ may be obtained
  as follows. For any fixed $\de>0$ we may find a finite
  open cover of $M$ by $\de$-balls:
  $\{B(x_i,\de),i=1,\dots,k\}$. We observe that since
  $\mu_0$ is a finite measure, there exist arbitrarily small
  values $\eta>0$ such that $\mu_0(\partial
  B(x_i,\de+\eta))=0$ for all $i=1,\dots,k$. Moreover we
  automatically have $m(\partial B(x_i,\de+\eta))=0$ also.
  Let us fix such a $\eta$.  Then the partition $\xi=
  \{B(x_1,\de+\eta),M\setminus B(x_1,\de+\eta)\} \vee
  \dots\vee \{B(x_k,\de+\eta), M\setminus B(x_k,\de+\eta)\}$ is as
  stated.
\end{proof}

\begin{theorem}
\label{thm.acumulacao}
In this setting, for all weak$^*$ accumulation point
$\mu_0$ of $\mu_t$ when $t\to 0^+$ we have
\begin{equation}
\label{eq.semicontentropy}
\limsup_{t\to 0^+} h_{\mu_t}(f_t) \le h_{\mu_0}(f_0).
\end{equation}
\end{theorem}

This result together with Ruelle's inequality will show that
every weak$^*$ accumulation point
$\mu_0$ of $\mu_t$ when $t\to 0^+$ satisfies the Entropy
Formula.

\begin{proof}
  Let us fix a weak$^*$ accumulation point $\mu_0$ of
  $\mu_t$ when $t\to 0^+$ and a partition $\xi$ as in
  Theorem~\ref{thm.particao}. Then by the Kolmogorov-Sinai
  Theorem~\cite{Man87} and setting $\xi_n=\vee_{j=0}^{n-1}
  f_t^{-j}\xi$ we have for any given fixed $n\ge1$
  \[
    h_{\mu_t}(f_t) =
    h_{\mu_t}(f_t,\xi)=\inf_{k\ge1}\, \frac1k
    H_{\mu_t}(\xi_k) \le
\frac1n\int-\log\mu_t(\xi_n(x))\,d\mu_t(x).
  \]
Now since the boundary of every element of $\xi$ has $\mu_0$
measure zero, 
then we have the following convergence
\[
\frac1n\int-\log\mu_t(\xi_n(x))\,d\mu_t(x) \to
\frac1n\int-\log\mu_0(\xi_n(x))\,d\mu_0(x) = \frac1n
H_{\mu_0}(\vee_{j=0}^{n-1} f_0^{-j}\xi).
\]
Since this holds for all $n\ge1$, we have
\[
\limsup_{t\to0^+} h_{\mu_t}(f_t) \le h_{\mu_0}(f_0),
\]
completing the proof.
\end{proof}

From Theorem~\ref{thm.acumulacao} we conclude that the
Entropy Formula~\eqref{eq:E-F} holds for every weak$^*$
accumulation point $\mu_0$ of $(\mu_t)_{t>0}$ when $t\to
0^+$, since as already observed the opposite inequality
in~\eqref{eq.semicontentropy} is always true by \cite{Ru}.

Finally, by Theorem~\ref{thm:main-FE}, we see that  every weak$^*$
accumulation point $\mu_0$ as above is the Dirac mass
$\de_0$, which ends the proof of Theorem~\ref{thm:mainstatistical}.



\section{Stochastic stability}
\label{sec:stochastic-stability}

Here we prove Theorem~\ref{thm:stochstabnosinks}.  We
consider the family $f_t(x)=f_0(x)+t$, where $f_0$ satisfies
(H1), which is a generic unfolding of the saddle-node at
$0$.  Hence for all $t>0$ close enough to $0$ the map $f_t$
is uniformly expanding, by Theorem~\ref{thm.expansao}.


\subsection{Uniqueness of stationary probability measures}
\label{sec:uniq-stat-prob}

We note that by the choice of the family
$(f_t)_{t\in[0,t_0]}$, generically unfolding the saddle-node
at $0$, we have that there exists $\zeta>0$ such that for
all $x\in M$
\begin{equation}
  \label{eq:zetabola}
  \{ f_t(x) : t\in\supp(\te_\ep)\} \supset B(f_{t^*}(x),\zeta),
\end{equation}
for some fixed $t^*\in\supp(\te_\ep)$, where $B(z,\zeta)$ is
the ball of radius $\zeta$ centered at $z$.
This holds just because $\supp(\te_\ep)$ has nonempty
interior and the map $t\mapsto f_t(x)$ is continuous (in
fact $C^2$) for every fixed $x$.

Let us define $f_x: [0,t_0]\to M, t\mapsto f_t(x)=f_0(x)+t$ for any
given fixed $x\in M$.  The condition $\te_\ep\ll m$ ensures
that for every $x\in M$ we have $(f_x)_* (\te_\ep^\bN) \ll
m$, where $(f_x)_* (\te_\ep)$ is the
probability measure defined by
\[
\int \vfi \, d(f_x)_* (\te_\ep) 
=  \int \vfi(f_{ t}(x)) \, d\te_\ep( t)
\]
for every bounded measurable function $\vfi:M\to\bR$.
Indeed, if $E$ is a Borel subset of $M$ such that $m(E)=0$, then 
\[
(f_x)_* (\te_\ep)(E)=\int  1_E(f_0(x)+t) \, d\te_\ep( t) = 
\int 1_{E-f_0(x)} \, d\te_\ep = \te_\ep(E-f_0(x))=0,
\]
because $m(E-f_0(x))=m(E)=0$. The definition of stationary
measure shows that every $\ep$-stationary measure $\mu^\ep$
is such that
\[
\int \vfi \, d\mu^\ep = \int \int \vfi(f_t(x))\, d\te_\ep(t)
\, d\mu^\ep(x) = \int [(f_x)_*\te_\ep]\vfi \,  d\mu^\ep(x),
\]
hence $\mu^\ep\ll m$ also. Since $\mu^\ep(\supp(\mu^\ep))=1$
we get $m(\supp(\mu^\ep))>0$ using the absolute continuity.

A standard property of $\ep$-stationary measures is that
$f_t(\supp(\mu^\ep))\subset\supp(\mu^\ep)$ for every
$t\in\supp(\te_\ep)$, see e.g.~\cite{Ar}.

This invariance property together with \eqref{eq:zetabola}
show that there exist $\zeta>0$ and $x\in\supp(\mu^\ep)$
such that $\supp(\mu^\ep)\supset B(f_{t^*}(x),\zeta)$.  Thus
the support of $\mu^\ep$ has nonempty interior.  By the
uniform expansion we know that $f_t$ is transitive (even
topologically mixing, see e.g.\cite{Man87}) for all
$t\in\supp(\te_\ep)$, hence we conclude that
\emph{$\supp(\mu^\ep)=M$ for every $\ep$-stationary
  probability measure $\mu^\ep$}.

Under the conditions assumed in the statement of
Theorem~\ref{thm:stochstabnosinks} together with
property~\eqref{eq:zetabola} it is known (see
e.g.~\cite{Ar}) that there are at most finitely many
$\ep$-stationary ergodic absolutely continuous probability
measures with pairwise disjoint supports. Since we have
shown that any $\ep$-stationary measure has full support in
$M$, we conclude that for every $\ep>0$ there is a unique
$\ep$-stationary absolutely continuous and ergodic measure
$\mu^\ep$, as stated in item (1) of
Theorem~\ref{thm:stochstabnosinks}.


\subsection{Entropy and random generating partitions}
\label{sec:entr-rand-gener}

Let $\mu^\ep$ be a $\ep$-stationary measure as defined 
above. Here we give two equivalent definitions of the
entropy of $\mu^\ep$ to be used in what follows.

\begin{theorem}{\cite[Thm. 1.3]{Ki86}} \label{thm.metr-entr-rand}
  For any finite measurable partition $\xi$ of $M$ the limit
 $$
 h_{\mu^\ep}( \xi) = \lim_{n \rightarrow \infty}
 \frac{1}{n} \int H_{\mu^\ep} \big(
 \bigvee_{k=0}^{n-1}(f^k_{\omega})^{-1} \xi \big) \, d
 \theta_\ep^\bN (\omega)
 $$
 exists.
\end{theorem}
This limit is called the \emph{entropy of the random
  dynamical system with respect to $\xi$ and to $\mu^\ep$}.
As in the deterministic case the above limit can be
replaced by the infimum.

The \emph{metric entropy} of the random
dynamical system is defined as
\[
h_{\mu^\ep} = \sup_\xi h_{\mu^\ep} ( \xi),
\]
where the supremum is taken over all finite measurable
partitions.

It seems natural to define the entropy of a random system by
$h_{\theta_{\ep}^\bN \times \mu^\ep} (S)$ where $S$ is the
skew-product map $S:[0,t_0]^\bN\times M\to[0,t_0]^\bN\times
M, \,\, (\omega,x)\mapsto (\sigma(\omega),f_{t_1}(x))$, and
$\sigma:[0,t_0]^\bN\to[0,t_0]^\bN$ is the left shift on
sequences.  However (see e.g.  \cite[Thm.  1.2]{Ki86}) under
some mild conditions the value of this function is infinite.
But the conditional entropy of $\theta_{\ep}^\bN \times
\mu^\ep$ with respect to a suitable $\sigma$-algebra of
subsets coincides with the entropy as defined above.

Let $\Omega=[0,t_0]^\bN$ be endowed with standard infinite
product (Tychonoff) topology, which makes $\Omega$ a compact
metric space. We consider the following compatible distance
in $\Omega$: given $\omega,\eta\in\Omega$
\[
D(\omega,\eta) = \sum_{i=1}^\infty
\frac1{2^i}\cdot d(\omega_i,\eta_i),
\]
where $\omega=(\omega_i)_{i\ge1}, \eta=(\eta_i)_{i\ge1}$ and
$d$ is the Euclidean distance on $[0,t_0]$.  Let $\cB$ be
the Borel $\sigma$-algebra of $\Omega$ and denote by $\cB
\times M$ the minimal $\sigma-$algebra containing all
products of the form $A \times M$ with $A \in \cB$.

In what follows we denote by $h_{
  \theta_{\ep}^\bN\times \mu^\ep}^{\cB \times M} (S)$ the
conditional metric entropy of $S$ with respect to the
$\sigma$-algebra $\cB \times M.$ (See e.g.~\cite{Bi65} for
definition and properties of conditional entropy.)

\begin{theorem}{\cite[Thm. 1.4]{Ki86}}
\label{thm.randentropyS}
Let $\mu^\ep$ be a $\ep$-stationary probability measure.
Then
\[
h_{\mu^\ep} = h_{ \theta_{\ep}^\bN
  \times \mu^\ep}^{ \cB \times M} (S).
\]
\end{theorem}

The Kolmogorov-Sinai result about generating partitions is
also available in a random version. We denote by
$\cA=\cB(M)$ the Borel $\sigma$-algebra of $M$ and say that
for a given fixed $\ep>0$, a finite partition $\xi$ is a
\emph{random generating partition for $\cA$} if
\begin{equation}
  \label{eq:randomgenerating}
  \bigvee_{k=0}^{+\infty} (f_{\omega}^k)^{-1} \xi = \cA \quad
  \text{for} \quad \theta_{\ep}^\bN-\text{almost all } \omega
  \in [0,t_0]^\bN.
  \end{equation}

\begin{theorem}{\cite[Cor. 1.2]{Ki86}}
  \label{thm.KSrandom}
  If $\xi$ is a random generating partition for $\cA$, then
  $h_{\mu^\ep}=h_{\mu^\ep}(\xi)$.
\end{theorem}


\subsubsection{Entropy Formula for random perturbations}
\label{sec:entr-form-rand}

We want to show that $\mu^\ep$ satisfies an Entropy
Formula analogous to~\eqref{eq:E-F} in the random setting.
The absolute continuity and ergodicity of $\mu^\ep$ gives
that $\mu^\ep$ satisfies the Entropy Formula in the
following form (see \cite{Li99}):
\begin{equation}
  \label{eq:34}
h_{\mu^\ep}=\lim_{n\to+\infty}\frac1n\log|
(f^n_\omega)'(x)|=
\lim_{n\to+\infty}\frac1n
\sum_{j=0}^{n-1} \log | f_0'(f_\omega^j(x)) |
=\int \log | f_0' | \, d\mu^\ep, 
  \end{equation}
  for $\te_\ep^\bN\times\mu^\ep$ almost every
  $(\omega,x)\in\Omega\times M$, as long as the random
  Lyapunov exponent given by the above limit is non-negative.
  (This limit does not depend on $(\omega,x)$ by the Ergodic
  Theorem.) Since by Theorem~\ref{thm.expansao} we have that
  the random Lyapunov exponent is non-negative for all
  $x\in\bS^1$ and $\theta_\epsilon^\bN$-a.e. $\omega$, then
  the Entropy Formula~\eqref{eq:34} holds.

\subsubsection{Constructing the generating partition}
\label{sec:constr-gener-part}

Here we use the previous results to prove the following
theorem analogous to Theorem~\ref{thm.particao}.

\begin{theorem}
\label{thm.particaoradom}
Let $\mu_0$ be a weak$^*$ accumulation point of $\mu^\ep$ when $\ep\to0$.
Then 
there exists a finite partition $\xi$ of $M$ which is
a $\mu^\ep\bmod0$
generating partition for all small enough $\ep>0$, and
also that $\mu_0(\partial\xi)=0$, i.e., the $\mu_0$ measure
of the boundary of the atoms of $\xi$ is zero.
\end{theorem}

\begin{proof}  
  A finite partition Lebesgue modulo zero whose boundary has
  also zero measure with respect to $\mu_0$ and with
  arbitrarily small diameter $\de>0$ may be obtained as
  already explained in the proof of
  Theorem~\ref{thm.particao}.

Now we show that if the diameter $\delta$ of $\xi$ satisfies
$0<\de<\de_1$, where $\de_1$ is the \emph{injectivity
  radius} of $f_t$ for all $t\in[0,t_0]$, i.e.,
\[
f_t\mid B(x,\delta_1)\quad\mbox{is a diffeomorphism onto its
  image,}\quad t\in[0,t_0], \, x\in\bS^1.
\]
(since $f_t$ is a family of local diffeomorphisms, this
value $\delta_1>0$ is guaranteed to exist), then $\xi$ is a random
generating partition for the Borel $\sigma$-algebra as in
\eqref{eq:randomgenerating} for all small enough $\ep>0$.

Indeed, let $x,y\in\bS^1$ be given  and let $\omega$ be a
$\theta_\epsilon^\bN$-generic sequence such that
\begin{equation}
  \label{eq:mais2}
  \dist \Big( f_\omega^n(x),f_\omega^n(y) \Big) \le \de
  \quad
\mbox{for every}\quad n\ge1,
\end{equation}
where $0<\delta<\delta_1$. Let $n_1<n_2<n_3<\dots$ be given
by Theorem~\ref{thm.expansao} and
Remark~\ref{rmk:mais1}. Then we have for all $n\ge1$
\[
\dist(x,y)\le \sigma_0^{-\sum_{k=1}^n (n_{2k}-n_{2k-1})}
\cdot \dist\Big( f_\omega^{n_{2k}}(x), f_\omega^{n_{2k}}(y)\Big)
\]
because, by assumption~\eqref{eq:mais2},
$f_\omega^{n_{2k}}(x), f_\omega^{n_{2k}}(y)$ are always in a
region where $f_t$ is invertible.

Hence for a partition $\xi$ with $\diam \xi<\delta_1$ and
$\mu_0(\partial\xi)=0$, setting
$\xi_{n,\omega}=\vee_{j=0}^{\,n-1}(f_\omega^j)^{-1}\xi$ we
have that for every $x\in\bS^1$
\[
\diam \xi_{n_{2k},\omega}(x) \to 0 \quad\mbox{when}\quad k\to\infty
\]
for $\theta_\epsilon^\bN$-a.e. $\omega$. This implies that
$\bigvee_{n\ge1}\xi_{n,\omega}=\cA, \, \mu^\ep\bmod0$,
finishing the proof.

\end{proof}


\subsection{Accumulation measures and Entropy Formula}
\label{sec:accum-meas-entr}

Now we prove that every weak$^*$ accumulation measure
$\mu_0$ of $(\mu^\ep)_{\ep>0}$ when $\ep\to0$ satisfies the
Entropy Formula. 

We start by fixing a weak$^*$ accumulation point $\mu_0$
of $\mu^\ep$ when $\ep\to 0$: there exists $\ep_k\to0$ when
$k\to\infty$ such that $\mu=\lim_k\mu^{\ep_k}$. We also fix
a uniform random generating partition $\xi$ as in the
previous subsection.

We need to construct a sequence of partitions of
$\Omega\times M$ according to the following result. 
We set $\omega_0=(0,0,0,\dots)\in\Omega$ in what follows.

\begin{lemma}
\label{lem.seqnpartitions} There exists an increasing
sequence of measurable partitions $(\cB_n)_{n\ge1}$ of
$\Omega$ such that
 \begin{enumerate}
 \item $\omega_0\in\inte (\cB_n(\omega_0))$ for all $n\ge1$;
 \item $\cB_n \nearrow \cB$, $\te^{\ep_k} \bmod 0$ for all
   $k\ge1$ when $n\to\infty$;
 \item $\lim_{n \to \infty} H_{\rho} (\xi \mid \cB_n)
 = H_{\rho} (\xi \mid \cB)$ for every measurable finite
 partition $\xi$ and any $S$-invariant probability measure $\rho$.
 \end{enumerate}  
\end{lemma}

\begin{proof}
  In this proof all distances and diameters are taken with
  respect to the distance $D$ on $\Omega$.

  For the first two items we let $\cC_n$ be a finite
$\te_{\ep_k}\bmod 0$ partition of $\Omega$ such that
$t_0\in\inte(\cC_n(t_0))$ with $\diam \cC_n\to0$ when
$n\to\infty$. Example: take a cover $(B(t,1/n))_{t\in X}$ of
$\Omega$ by $1/n$-balls and take a subcover $U_1,\dots,U_k$ of
$\Omega\setminus B(t_0,2/n)$ together with $U_0=B(t_0,3/n)$; then
let $\cC_n=\{U_0,M\setminus U_0\}\vee\dots\vee \{U_k,
M\setminus U_k\}$.

We observe that we may assume that the boundary of these
balls has null $\te_{\ep_k}$-measure for all $k\ge1$, since
$(\te_{\ep_k})_{k\ge1}$ is a denumerable family of
non-atomic probability measures on $\Omega$.  Now we set
\[
\cB_n=\cC_n\times\stackrel{n}{\dots}\times\cC_n\times\Omega\quad
\mbox{for all  }n\ge1.
\]
Then since $\diam \cC_n\le 2/n$ for all $n\ge1$ we have that
$\diam \cB_n\le 2/n$ also and so tends to zero when
$n\to\infty$.  Clearly $\cB_n$ is an increasing sequence of
partitions. Hence $\vee_{n\ge1} \cB_n$ generates the
$\sigma$-algebra $\cB$, $\te^{\ep_k} \bmod0$ (see e.g.
\cite[Lemma 3, Chpt. 2]{Bi65}) for all $k\ge1$. This proves
items (1) and (2).

Item (3) of the statement of the lemma is Theorem 12.1 of
Billingsley~\cite{Bi65}.
\end{proof}

Now we use some properties of conditional entropy to obtain
the right inequalities. We start with
\begin{eqnarray*}
  h_{\mu^{\ep_k}}
  &=&
  h_{\mu^{\ep_k}}(\xi) = h_{
  \theta_{{\ep_k}}^\bN \times \mu^{\ep_k}}^{ \cB \times M} (S,
  \Omega\times\xi)
  \\
  &=&
  \inf_{n\ge1}\frac1n H_{\te_{\ep_k}^\bN\times\mu^{\ep_k}}\left(
  \bigvee_{j=0}^{n-1} (S^j)^{-1}(\Omega\times\xi) \mid
  \cB\times M \right)
\end{eqnarray*}
where the first equality comes from the random
Kolmogorov-Sinai Theorem~\ref{thm.KSrandom} and the second
one can be found in Kifer~\cite[Thm. 1.4, Chpt. II]{Ki86},
with $\Omega\times\xi=\{ \Omega\times A: A\in\xi\}$.  Hence
for arbitrary fixed $N\ge1$ and for any $m\ge1$
\begin{eqnarray*}
 h_{\mu^{\ep_k}}
 &\le&
 \frac1N \cdot H_{\te_{\ep_k}^\bN\times\mu^{\ep_k}}\left(
  \bigvee_{j=0}^{N-1} (S^j)^{-1}(\Omega\times\xi) \mid
  \cB\times M \right)
 \\
 &\le&
 \frac1N \cdot H_{\te_{\ep_k}^\bN\times\mu^{\ep_k}}\left(
  \bigvee_{j=0}^{N-1} (S^j)^{-1}(\Omega\times\xi) \mid
  \cB_m\times M \right)
\end{eqnarray*}
because $\cB_m\times M\subset\cB\times M$. Now we fix $N$ and $m$,
let $k\to \infty$ and note that since
$\mu_0(\partial\xi)=0=\de_{\omega_0}(\partial \cB_m)$ it must be
that
\[
(\de_{\omega_0}\times\mu_0)(\partial(B_i\times\xi_j))=0\quad
\mbox{for all  } B_i\in\cB_m \mbox{  and  } \xi_j\in\xi,
\]
where $\de_{\omega_0}$ is the Dirac mass concentrated at
$\omega_0\in\Omega$.  Thus we get by weak$^*$ convergence of
$\te_{\ep_k}^\bN\times\mu^{\ep_k}$ to
$\delta_{\omega_0}\times\mu_0$ when $k\to\infty$
\begin{equation}
  \label{eq:5a}
  \limsup_{k\to\infty} h_{\mu^{\ep_k}}
  \le
  \frac1N \cdot H_{\de_{\omega_0}\times\mu_0}\left(
  \bigvee_{j=0}^{N-1} (S^j)^{-1}(\Omega\times\xi) \mid
  \cB_m\times M \right)
  =\frac1N \cdot H_{\mu_0} \big(\bigvee_{j=0}^{N-1} f^{-j}\xi \big).
\end{equation}
Here it is easy to see that the middle conditional entropy
of~\eqref{eq:5a} (involving only finite partitions) equals
\[
\frac1N \sum_i \mu_0(P_i)\log\mu_0(P_i),
\]
with $P_i=\xi_{i_0}\cap f^{-1}\xi_{i_1}\cap\dots\cap
f^{-(N-1)} \xi_{i_{N-1}}$ ranging over all possible
sequences $\xi_{i_0},\dots,\xi_{i_{N-1}}\in \xi$.
Finally, since $N$ was an arbitrary integer,
it follows from~\eqref{eq:34}, \eqref{eq:5a} and the Ruelle
Inequality that
\[
\int\log|f_0'|\,d\mu_0
\le\limsup_{k\to\infty} h_{\mu^{\ep_k}}\le h_{\mu_0}(f_0)
\le \int\log|f_0'|\,d\mu_0,
\]
showing that $\mu_0$ satisfies the Entropy Formula.

To conclude the proof of Theorem~\ref{thm:stochstabnosinks}
we observe that $\mu_0$ is $f_0$-invariant by construction
and since it satisfies the Entropy Formula, Theorem~\ref{thm:main-FE}
ensures that $\mu_0=\de_0$ the Dirac mass at the saddle-node
$0$.


\bibliographystyle{abbrv}

\end{document}